\documentclass{article}
\usepackage{amsmath,amssymb,amsthm,latexsym}

\newtheorem{theorem}{Theorem}
\newtheorem{lemma}{Lemma}

\begin{document}
\title{The Lattice Dimension of a Tree}
	\author{Sergei~Ovchinnikov \\
	Mathematics Department\\
	San Francisco State University\\
	San Francisco, CA 94132\\
	sergei@sfsu.edu} 
\date{\empty}
\maketitle

\begin{abstract}
The lattice dimension of a graph $G$ is the minimal dimension of a cubic lattice in which $G$ can be isometrically embedded. In this note we prove that the lattice dimension of a tree with $n$ leaves is $\lceil n\slash 2\rceil$.
\end{abstract}

In the paper, $T$ is a tree with $q$ edges and $n$ leaves. By the Djokovi\'{c} theorem~\cite{dD73} (see also~\cite{aA98,mD97,wI00}), $T$ can be isometrically embedded into the $q$--cube $\mathcal{Q}_q$ and cannot be embedded into a cube of a smaller dimension. The number $q$ is the \emph{isometric dimension} $\text{dim}_I(T)$ of the tree $T$~\cite{mD97}. Thus the vertices of $T$ can be labeled with $0\slash 1$ $q$--dimensional vectors in such a way that the distance between two edges is the Hamming distance between the corresponding vectors. This labeling is also known as an \emph{addressing scheme} for $T$~\cite{aA98}.

Let $\mathbb{Z}^d$ be the edge skeleton of the standard cubic lattice of dimension $d$. Simple examples show that, generally speaking, $T$ can be isometrically embedded into $\mathbb{Z}^d$ of a smaller dimension than $\text{dim}_I(T)$. In a more general setting, one can consider the following problem: for a given graph $G$, find the minimum possible dimension $d$ such that $G$ can be isometrically embedded into $\mathbb{Z}^d$. We call this number $d$ the \emph{lattice dimension} of $G$ and denote it by $\text{dim}_Z(G)$. Recently, Eppstein~\cite{dE04} showed how to express the lattice dimension of a partial cube $G$ in terms of a maximal matching in the semicube graph $\text{Sc}(G)$ associated with $G$, and presented a polynomial time algorithm for finding this dimension. In this paper, we prove directly that the lattice dimension of a tree with $n$ leaves is $\lceil n\slash 2\rceil$.

Clearly, a tree with two leaves is embeddable into $\mathbb{Z}$. In what follows $T$ is a tree with $n>2$ leaves. For a given leaf $v\in T$ we denote $v_T$ the closest vertex in $T$ which has a degree greater than two and call the unique path $(v_T,v)$ the \emph{hanging path} to $v$.

\begin{lemma} \label{StarEmbedding}
Suppose the star $K_{1,n}$ is isometrically embedded into $\mathbb{Z}^d$. Then
$$
d\geq \lceil n\slash2\rceil,
$$
Moreover, $\text{\emph{dim}}_Z(K_{1,n})=\lceil n\slash2\rceil$.
\end{lemma}

\begin{proof}
We may assume that the center of $K_{1,n}$ is represented by the zero vertex in $\mathbb{Z}^d$. Then each leaf in $K_{1,n}$ is represented by the end vertex of a coordinate unit vector $e_k$ or $-e_k$. It follows that $2d\geq n$. Since $d$ is an integer, we have $d\geq \lceil n\slash2\rceil$.

We construct an isometric embedding of $K_{1,n}$ into $\mathbb{Z}^{\lceil n\slash2\rceil}$ as follows. Let $v_1,\ldots,v_n$ be the leaves and $c$ be the center of $K_{1,n}$. We map $c$ into the zero vertex of $\mathbb{Z}^{\lceil n\slash2\rceil}$, $v_1$ and $v_2$ into the end points of $e_1$ and $-e_1$, respectively, and so on. Clearly, we obtain an isometric embedding of $K_{1,n}$ into $\mathbb{Z}^{\lceil n\slash2\rceil}$.
\end{proof}

\begin{lemma} \label{TreeEmbedding}
Let $T$ be a tree with $n$ leaves. Suppose that $T$ is isometrically embedded into $\mathbb{Z}^d$. Then
$$
d\geq \lceil n\slash2\rceil.
$$
\end{lemma}

\begin{proof}
Let $(u,v)$ be an inner edge of $T$, that is an edge with vertices that are not leaves. Vertices $u$ and $v$ are represented by vertices in $\mathbb{Z}^d$ that are different only in one position, say,
$$
u=(x_1,x_2,\ldots,x_d),\quad v=(x'_1,x_2,\ldots,x_d),
$$
with $|x_1-x'_1|=1$.
Let us select all edges in $\mathbb{Z}^d$ with vertices that have the same first coordinates as $u$ and $v$. Since $T$ is a tree and the embedding is isometric, the edge $(u,v)$ is the only edge among selected which belongs to $T$. Now we remove all the selected edges in $\mathbb{Z}^d$ and, for each selected edge, identify vertices defining this edge. We obtained a tree which is a contraction~\cite{wI00} of $T$ having the same number $n$ of leaves and embedded into another copy of $\mathbb{Z}^d$. By repeating this process, we end up with the star $K_{1,n}$ embedded into $\mathbb{Z}^d$. By Lemma~\ref{StarEmbedding}, we have $d\geq\lceil n\slash 2\rceil$.
\end{proof}

Now we prove the main theorem of the paper.
\begin{theorem} \label{TreeDimension}
The lattice dimension of a tree $T$ with $n$ leaves is $\lceil n\slash 2\rceil$.
\end{theorem}

\begin{proof}
By Lemma~\ref{TreeEmbedding}, it suffices to show that there is an isometric embedding of $T$ into $\mathbb{Z}^{\lceil n\slash 2\rceil}$. To construct this embedding we use an inductive argument.

Suppose that the statement of the theorem is true for all trees with the number of leaves less than $n$ for a given $n>2$. Let $v$ and $u$ be two leaves of $T$ and $(v_T,v)$ and $(u_T,u)$ be their hanging paths. The degrees of vertices $v_T$ and $u_T$ are greater than two and these two vertices are not necessarily different. 

Suppose first that $v_T\not=u_T$. By removing the hanging paths $(u_T,u)$ and $(v_T,v)$, we obtain a new tree $T'$ with $(n-2)$ leaves. By the induction hypothesis, $T'$ can be isometrically embedded into $\mathbb{Z}^{d}$ with $d=\lceil\frac{n-2}{2}\rceil=\lceil \frac{n}{2}\rceil-1$. Now we embed, in a natural way, $\mathbb{Z}^d$ into $\mathbb{Z}^{d+1}$ and add the paths $(u_T,u)$ and $(v_T,v)$ to the image of $T'$ in the positive and negative directions of the new dimension, respectively. Clearly, we constructed an isometric embedding of $T$ into $\mathbb{Z}^{\lceil n\slash 2\rceil}$.

Suppose now that all leaves $v$ have the same vertex $c=v_T$ in paths $(v_T,v)$. Then $c$ is the only vertex in $T$ of degree greater than two. By modifying in an obvious way the construction from Lemma~\ref{StarEmbedding}, we obtain the required embedding.
\end{proof}

Suppose $T$ is isometrically embedded into $\mathbb{Z}^d$ with $d=\lceil n\slash 2\rceil$. The projection of $T$ into the $k$'th coordinate axis is a path $P_{l_k}$ of length $l_k>0$. Thus $T$ is actually embedded into the \emph{complete grid graph}~\cite{wI00}
\begin{equation*}
P=P_{l_1}\Box P_{l_2}\Box\cdots\Box P_{l_d},
\end{equation*}
where $\Box$ denotes the Cartesian product. We have
\begin{equation*}
d = \text{dim}_Z(T)=\text{dim}_Z(P)\quad\text{and}\quad q = \text{dim}_I(T)=\text{dim}_I(P).
\end{equation*}

Although the grid $P$ is not uniquely defined by $T$, the set $\mathcal{L}=\{l_1,\cdots,l_d\}$ determines both the lattice and isometric dimensions of $T$. Namely, 
\begin{equation*}
\text{dim}_Z(T)=|\mathcal{L}|\quad\text{and}\quad \text{dim}_I(T)=l_1+\cdots+l_d.
\end{equation*}

Clearly, we may assume that $P_{l_k}=\{0,1,\cdots,l_k\}$. Then an isometric embedding of $T$ into $P$ defines an addressing scheme for $T$ in which each vertex is labeled by a $d$--dimensional vector with the $k$'s coordinate ranging from $0$ to $l_k$. The distance between two edges is the $\ell_1$-distance between the corresponding vectors.

\end{document}